\documentclass[11pt,a4paper]{amsart}
\usepackage[T1]{fontenc}
\usepackage{lmodern}
\usepackage{amsmath,amssymb,mathtools}
\usepackage[margin=28mm]{geometry}
\usepackage{microtype}
\usepackage{booktabs}
\usepackage{enumitem}
\usepackage{needspace}
\usepackage[hidelinks]{hyperref}
\hypersetup{
  pdftitle={A sharp weighted L2 Bernstein inequality for partial derivatives on the unit ball},
  pdfauthor={Tomasz Beberok},
  pdfsubject={Polynomial inequalities, exact constants and extremal polynomials}
}
\numberwithin{equation}{section}
\newtheorem{theorem}{Theorem}[section]
\newtheorem{proposition}[theorem]{Proposition}
\newtheorem{corollary}[theorem]{Corollary}
\newtheorem{lemma}[theorem]{Lemma}
\theoremstyle{definition}

\theoremstyle{remark}
\newtheorem{remark}[theorem]{Remark}
\newcommand{\R}{\mathbb R}

\newcommand{\B}{\mathbb B}
\newcommand{\Sph}{\mathbb S}
\newcommand{\HH}{\mathcal H}
\newcommand{\PP}{\Pi}
\newcommand{\dd}{\,\mathrm d}
\newcommand{\ip}[2]{\langle #1,#2\rangle_\mu}
\newcommand{\nm}[1]{\lVert #1\rVert_\mu}

\newcommand{\doi}[1]{\href{https://doi.org/#1}{\nolinkurl{#1}}}
\setlist[enumerate]{label=\textup{(\roman*)},leftmargin=2em,itemsep=3pt}
\allowdisplaybreaks[1]
\title[Sharp Bernstein inequality for partial derivatives]
{A sharp weighted $L^2$ Bernstein inequality for partial derivatives on the unit ball}
\author{Tomasz Beberok}
\address{Department of Applied Mathematics,
University of Agriculture in Krakow,
Krakow, Poland}
\email{tomasz.beberok@urk.edu.pl}
\date{}
\subjclass[2020]{41A17, 33C45, 42C10}
\keywords{Bernstein inequality, unit ball, Jacobi polynomials, sharp constant, extremal polynomials}

\begin{document}
\begin{abstract}
We determine the sharp constant in the weighted $L^2$
Bernstein inequality for coordinate derivatives on the
unit ball with a Jacobi weight and the strengthened
derivative factors introduced by Beberok and Xu.
The constant is expressed in terms of the largest zeros
of Jacobi polynomials, and all extremal polynomials are
characterized explicitly. The proof uses an orthogonal
polynomial decomposition that gives diagonal formulas
for both the squared norm and the weighted derivative
integral.
We also obtain two-sided bounds for the sharp constant
and determine its leading asymptotic behavior as the
polynomial degree tends to infinity.
\end{abstract}
\maketitle

\section{Introduction and the sharp inequality}\label{sec:intro}

For $a,b > -1$ and $n \ge 1$, the classical weighted $L^2$
Bernstein inequality for Jacobi weights states that
\[
 \int_{-1}^{1} (1-t)^{a+1} (1+t)^{b+1} |p'(t)|^2 \,\dd t
 \le n(n+a+b+1)
      \int_{-1}^{1} (1-t)^a (1+t)^b |p(t)|^2 \,\dd t
\]
for every real polynomial $p$ of degree at most $n$.
The constant is sharp, and equality holds for the Jacobi polynomial $P_n^{(a,b)}$; see \cite[Theorem~2.1]{GM}.
In several variables, weighted Bernstein inequalities
on the sphere and the ball were studied by Dai \cite{Dai2006}
for doubling weights.
Dai and Xu \cite{DX2010} developed moduli of smoothness
based on angular and weighted coordinate derivatives
and established direct and inverse theorems for
polynomial approximation on these domains.
These results explain the role of the derivative weights
in measuring smoothness and estimating approximation errors.

Let $\B^d=\{x \in \R^d: \|x\| \le 1\}$,
 $W_\mu(x)=(1-\|x\|^2)^\mu$,
 $\mu>-1$,
where $\|x\|$ denotes the Euclidean norm.
We write $\PP_n^d$ for the real polynomials in $d$ variables of
total degree at most $n$, and set
\begin{equation}\label{eq:inner}
 \ip{f}{g} = \int_{\B^d} f(x) g(x )W_\mu(x) \dd x,\qquad
 \nm{f}^2 = \ip{f}{f}.
\end{equation}
Throughout the paper, $d \ge 2$ and $n \ge 1$, unless stated otherwise.
For $1\le i\le d$, we write $\partial_i p$ for the
partial derivative of $p$ with respect to the coordinate $x_i$.
Bernstein inequalities on the ball have been studied for both coordinate
and angular derivatives.  For the Jacobi weight, the classical weighted Bernstein inequality on the ball gives
\begin{equation}\label{eq:classicalball}
 \sum_{i=1}^d \nm{\sqrt{1-\|x\|^2}\,\partial_i p}^{2}
 + \sum_{1\le i<j\le d}  \nm{D_{ij}p}^{2}
 \le n(n+2\mu+d)         \nm{p}^{2},
\end{equation}
where $D_{ij} = x_i \partial_j - x_j \partial_i$; see
\cite[Chapter~11]{DX} and \cite[Theorem~4.1]{BX}.
Kro\'o \cite[Theorem~1]{Kroo} determined the sharp constant
for the first sum alone and showed that its value depends
on the parity of $n$.
The corresponding extremal polynomials were characterized
in \cite{Kroo2}.
Beberok and Xu \cite[Theorem~3.1]{BX} established stronger coordinate inequalities involving the factors
\begin{equation}\label{eq:Phi}
 \Phi_i(x)=
 \frac{\sqrt{1-\|x\|^2}}{\sqrt{1-\|x\|^2+x_i^2}},
 \qquad \|x\| < 1,\quad 1 \le i \le d.
\end{equation}
Their result implies, in particular,
\begin{equation}\label{eq:BXbound}
 \nm{\Phi_i\partial_i p} \le c n \nm{p},
 \qquad p \in \PP_n^d,
\end{equation}
with $c$ independent of $n$.
Since $\Phi_i(x) \ge \sqrt{1- \|x\|^2}$, this controls a stronger weighted
partial derivative than the corresponding individual term in \eqref{eq:classicalball}.
Stronger derivative weights have been obtained for several multivariate domains.
Xu \cite{Xu2023} established strengthened Bernstein inequalities on conic domains and triangles.
Ge and Xu \cite{GX} obtained strengthened inequalities on the simplex for doubling weights.
On a trapezoidal domain, Beberok \cite{TB} obtained sharp weighted $L^2$ Bernstein inequalities
for directional derivatives and strengthened estimates for coordinate derivatives.

The purpose of the present work is to determine the smallest constant
$C_{n,d,\mu}$ for which
\begin{equation}\label{eq:desired}
 \nm{\Phi_i \partial_i p}
                          \le C_{n,d,\mu} \nm{p},
 \qquad p \in \PP_n^d,
\end{equation}
and to identify all polynomials attaining equality.
In other words, we seek to evaluate
\begin{equation}\label{eq:Cdef}
 C_{n,d,\mu}^2=
 \max_{0 \ne p \in \PP_n^d}
 \frac{\displaystyle \int_{\B^d}
       \frac{1-\|x\|^2}{1-\|x\|^2 + x_i^2}
       |\partial_i p(x)|^2 W_\mu(x) \dd x}
      {\displaystyle\int_{\B^d} |p(x)|^2 W_\mu(x) \dd x}.
\end{equation}
The value is independent of $i$ by invariance under coordinate permutations.
We will reduce this extremal problem to a finite family of one-dimensional inequalities.

For each integer $k\ge1$, let
\[
 \mathcal L_k=
 \begin{cases}
  \{0,1\},&k=1,\\
  \{0,1,2,\ldots\},&k\ge2.
 \end{cases}
\]
To state the sharp constant, define the index set
\[
 \mathcal I_{n,d}
 =\{(j,\ell): 1 \le j \le n,\, \ell \in \mathcal L_{d-1},\,
                  j + \ell \le n\}.
\]
For $(j,\ell)\in\mathcal I_{n,d}$, let
\begin{equation}\label{eq:parameters}
 m_{j,\ell} = \left\lfloor\frac{n-j-\ell}{2}\right\rfloor, \quad
 \alpha_j   = \mu + j - \frac12, \quad
 \beta_\ell = \ell + \frac{d-3}{2}, \quad
 \eta_j     = j(j+2\mu+1).
\end{equation}
Let $N_{j,\ell} = m_{j,\ell} + 1$ and let
$
 -1< \xi_{j,\ell,1} <\cdots<
 \xi_{j,\ell,N_{j,\ell}} <1
$
be the zeros of $P_{N_{j,\ell}}^{(\alpha_j,\beta_\ell)}$.
Here $P_k^{(\alpha,\beta)}$ denotes the classical Jacobi polynomial of degree $k$.
Since $\alpha_j,\beta_\ell > -1$ for every
$(j,\ell) \in \mathcal I_{n,d}$, the polynomial
$P_{N_{j,\ell}}^{(\alpha_j,\beta_\ell)}$
has exactly $N_{j,\ell}$ simple zeros in $(-1,1)$
by \cite[Theorem~3.3.1]{Szego}.

\begin{theorem}\label{thm:main}
Let $d\ge2$ and $n \ge 1$ be integers, and let $\mu > -1$.
The constant $C_{n,d,\mu}$ defined by
\eqref{eq:Cdef} satisfies
\begin{equation}\label{eq:sharpconstant}
 C_{n,d,\mu}^2=
 \max_{(j,\ell) \in \mathcal I_{n,d}}
 \frac{2j(j+2\mu+1)}
      {1 - \xi_{j,\ell,N_{j,\ell}}}.
\end{equation}
\end{theorem}

The nonzero polynomials attaining equality in
\eqref{eq:desired} are characterized in
Corollary~\ref{cor:extremizers}.

The connection between sharp polynomial inequalities and zeros of
orthogonal polynomials is classical; see, for example, Nikolov and Shadrin \cite{NS}.
For the factor \eqref{eq:Phi}, an expansion in Jacobi polynomials along
coordinate sections separates the differentiated variable.
A harmonic decomposition in the other $d-1$ variables then reduces
the quotient in \eqref{eq:Cdef} to quotients of two Jacobi integrals whose weights differ by the factor $1-u$.
Gaussian quadrature gives their exact extrema and the corresponding polynomials.

The one-dimensional facts and polynomial decompositions needed in the proof are established in Section~\ref{sec:lemmas}.
Section~\ref{sec:sharp} proves Theorem~\ref{thm:main} and determines all equality cases.
At the end of Section~\ref{sec:sharp}, we prove that
\begin{equation}\label{eq:introbounds}
 n(n+2\mu+1) \frac{n + \mu + d/2}{n + \mu + 1/2}
 \le C_{n,d,\mu}^2 \le n(n + 2\mu + d).
\end{equation}
Consequently $C_{n,d,\mu}/n \to 1$ for fixed $d,\mu$.

\section{Auxiliary lemmas}
\label{sec:lemmas}

Let $\{q_j\}_{j \ge 0}$ be the orthonormal Jacobi polynomials for the weight $(1 - t^2)^\mu$ on $[-1,1]$, with positive leading coefficients.  Equivalently,
\[
 q_j(t)=
 \frac{P_j^{(\mu,\mu)}(t)}
 {\left(\int_{-1}^1
       \big[P_j^{(\mu,\mu)}(v)\big]^2 (1-v^2)^\mu \dd v \right)^{1/2}}.
\]
In particular,
\begin{equation}\label{eq:qorth}
 \int_{-1}^1 q_j(t) q_k(t) (1-t^2)^\mu \dd t = \delta_{jk}.
\end{equation}
The following lemma gives the corresponding orthogonality relation
for the derivatives of $q_j$ with respect to the weight
$(1-t^2)^{\mu+1}$.
\begin{lemma}\label{lem:jacobi}
For $j,k \ge 0$,
\begin{equation}\label{eq:qenergy}
 \int_{-1}^1 q_j'(t) q_k'(t) (1-t^2)^{\mu+1} \dd t
 =\eta_j \delta_{jk},\qquad
 \eta_j = j(j + 2\mu + 1).
\end{equation}
Consequently, if $h = \sum\limits_{j=0}^m a_j q_j$, then
\begin{align}
 \int_{-1}^1 (h'(t))^2 (1-t^2)^{\mu+1} \dd t
   =\sum_{j=1}^m \eta_j a_j^2.\label{eq:oneDenergy}
\end{align}
For $m \ge 1$, this gives the sharp Bernstein inequality with
constant squared $\eta_m$, and its nonzero equality polynomials
are precisely the scalar multiples of $q_m$.
\end{lemma}

\begin{proof}
Since $q_j$ is a constant multiple of $P_j^{(\mu,\mu)}$,
the Jacobi differential equation
\cite[Theorem~4.2.1]{Szego} yields
\[
 (1-t^2) q_j''(t) - 2(\mu + 1)t q_j'(t) + \eta_j q_j(t)=0.
\]
Multiplying the above equation by $(1-t^2)^\mu$, we obtain
\[
 -\frac{\dd}{\dd t}
 \bigl((1-t^2)^{\mu+1} q_j'(t) \bigr)
 =\eta_j(1-t^2)^\mu q_j(t).
\]
Multiplying both sides by $q_k(t)$
and integrating by parts over $[-1,1]$,
we obtain
\begin{align*}
 \int_{-1}^{1}
      q_j'(t) q_k'(t) (1-t^2)^{\mu+1} \dd t
 =\eta_j \int_{-1}^{1}
      q_j(t) q_k(t) (1-t^2)^\mu\dd t.
\end{align*}
Hence \eqref{eq:qorth} yields \eqref{eq:qenergy}. The identity \eqref{eq:oneDenergy}
follows from  \eqref{eq:qenergy}.
Finally, for $j \ge 0$,
$
 \eta_{j+1} - \eta_j = 2(j+\mu+1) > 0.
$
Thus
\[
\sum_{j=1}^m \eta_j a_j^2 \le \eta_m \sum_{j=0}^m a_j^2,
\]
and equality for $m \ge 1$ holds exactly when $a_0 = \cdots = a_{m-1}=0.$
\end{proof}

The following lemma is a consequence of the classical Gaussian
quadrature formula; see \cite{Gautschi}.
We include a proof that also identifies all cases of equality.

\begin{lemma}\label{lem:gauss}
Let $\alpha,\beta > -1$, $m \ge 0$ be an integer, and
\[
 \omega(u)= u^\beta (1-u)^\alpha,\qquad
 \mathrm{P}(u) = P_{m+1}^{(\alpha,\beta)}(2u-1).
\]
Let $0 < u_1 < \cdots < u_{m+1} < 1$ be the zeros of $\mathrm{P}$.
For $1 \le \tau \le m + 1$, define
\begin{equation}\label{eq:cardinal}
 L_\tau(u)=\frac{\mathrm{P}(u)}{(u-u_\tau)\mathrm{P}'(u_\tau)},\qquad
 \gamma_\tau = \int_0^1 (L_\tau(u))^2 \omega(u) \dd u.
\end{equation}
Then, for all $h,q \in \PP_m^1$,
\begin{align}
 &\int_0^1 h(u) q(u) \omega(u)\dd u
   =\sum_{\tau=1}^{m+1} \gamma_\tau h(u_\tau) q(u_\tau),
   \label{eq:gaussfirst}\\
 &\int_0^1 (1-u) h(u) q(u) \omega(u) \dd u
   =\sum_{\tau=1}^{m+1} \gamma_\tau (1-u_\tau) h(u_\tau) q(u_\tau).
   \label{eq:gausssecond}
\end{align}
In particular, for every $h \in \PP_m^1$,
\begin{equation}\label{eq:gaussineq}
 \int_0^1 (h(u))^2 \omega(u)\dd u
 \le \frac1{1-u_{m+1}}
       \int_0^1 (1-u) (h(u))^2 \omega(u) \dd u.
\end{equation}
The constant is sharp, and a nonzero $h$ attains equality if and only if it is a scalar multiple of $L_{m+1}$.
\end{lemma}

\begin{proof}
The factor $u - u_\tau$ divides $\mathrm{P}$.
Hence $L_\tau$ is a polynomial of degree $m$, and
\[
 L_\tau(u_r)= \delta_{\tau r} \qquad (1 \le \tau,r \le m + 1).
\]
In particular, the polynomials $L_1,\ldots,L_{m+1}$
form a basis of $\PP_m^1$. Since $L_\tau(u_\tau)=1$, $L_\tau^2 - L_\tau$ vanishes at every zero of $\mathrm{P}$.
Therefore, for $m \ge 1$, there exists a polynomial $Q_\tau \in \PP_{m-1}^1$ such that
$
 L_\tau^2 - L_\tau = \mathrm{P} Q_\tau.
$
By the orthogonality of $\mathrm{P}$,
\[
 \int_0^1 \left( (L_\tau(u))^2 - L_\tau(u) \right) \omega(u) \dd u
 =\int_0^1 \mathrm{P}(u) Q_\tau(u) \omega(u) \dd u=0.
\]
For $m=0$, the same identity follows from $L_1 = 1$.
Thus
\begin{equation}\label{eq:gammaidentity}
 \gamma_\tau = \int_0^1 L_\tau(u) \omega(u) \dd u.
\end{equation}

Let $G \in \PP_{2m+1}^1$ and define
\[
 B(u) = \sum_{\tau=1}^{m+1} G(u_\tau) L_\tau(u).
\]
Then $B \in \PP_m^1$ and $B(u_\tau) = G(u_\tau)$ for every $\tau$.
Hence, $G-B$ is divisible by $\mathrm{P}$.
Consequently, there exists $A \in \PP_m^1$ such that
$
 G = A \mathrm{P} + B.
$
Since $A \in \PP_m^1$,
\[
 \int_0^1 A(u) \mathrm{P}(u) \omega(u) \dd u=0.
\]
Therefore, by \eqref{eq:gammaidentity} and the definition of $B$,
\begin{equation}\label{eq:gaussrule}
 \int_0^1 G(u) \omega(u) \dd u
   = \int_0^1 B(u) \omega(u) \dd u
   = \sum_{\tau=1}^{m+1} \gamma_\tau G(u_\tau).
\end{equation}
Applying \eqref{eq:gaussrule} to the polynomials $hq$ and
$(1-u)hq$, whose degrees are at most $2m$ and $2m+1$,
respectively, yields \eqref{eq:gaussfirst}
and \eqref{eq:gausssecond}.
In particular,
\begin{align}
 &\int_0^1 L_\tau(u) L_r(u) \omega(u) \dd u
   =\gamma_\tau \delta_{\tau r}, \notag\\
 &\int_0^1 (1-u) L_\tau(u) L_r(u) \omega(u) \dd u
   =\gamma_\tau (1 - u_\tau) \delta_{\tau r}.
   \label{eq:cardinalorthogonality}
\end{align}

Finally, let $h \in \PP_m^1$. Applying \eqref{eq:gaussfirst}
and \eqref{eq:gausssecond} with $q=h$, we obtain
\[
\begin{aligned}
\int_0^1  \frac{1-u}{1-u_{m+1}} (h(u))^2 \omega(u) \dd u
   -\int_0^1 (h(u))^2 \omega(u) \dd u
 &=
 \sum_{\tau=1}^{m+1} \gamma_\tau
 \left(\frac{1-u_\tau}{1-u_{m+1}}-1\right) (h(u_\tau))^2\\
 &=
 \sum_{\tau=1}^{m+1} \gamma_\tau
 \frac{u_{m+1}-u_\tau}{1-u_{m+1}} (h(u_\tau))^2.
\end{aligned}
\]
This proves \eqref{eq:gaussineq}.
All coefficients in the above sum except the last one are strictly positive.
Equality is therefore equivalent to $h(u_1) = \cdots = h(u_m)=0$.
If $h \in \PP_m^1$ and $h(u_r)=0$ for $1 \le r \le m$,
then the Lagrange interpolation formula gives
\[
 h(u) = \sum_{r=1}^{m+1} h(u_r) L_r(u)
      = h(u_{m+1}) L_{m+1}(u).
\]
For $m=0$, we have $L_1=1$, and the equality statement
holds for every nonzero constant polynomial.
\end{proof}

Let
 $
  x=(s,y) \in \R \times \R^{d-1}$,
 $
  a(y)=\sqrt{1-\|y\|^2}.
 $
For fixed $\|y\| < 1$, the coordinate section of $\B^d$ is $-a(y) \le s \le a(y)$.
Recall that
\[
 \Phi_i(x)=
 \sqrt{\frac{1 - \|x\|^2}{1 - \|x\|^2 + x_i^2}}.
\]
It is enough to consider $i=1$.  For $x=(s,y)$, this gives
\[
 (\Phi_1(s,y))^2
 =\frac{1 - s^2 - \|y\|^2}{1 - \|y\|^2}
 =1-\frac{s^2}{(a(y))^2}.
\]
The substitution $s=a(y)t$ therefore yields
\begin{equation}\label{eq:sectionchange}
 1 - s^2 - \|y\|^2 = (a(y))^2 (1-t^2),\qquad
 (\Phi_1(s,y))^2 = 1 - t^2,\qquad
 \dd s = a(y) \dd t.
\end{equation}
In particular,
$
 W_\mu(s,y) \dd s
 =(a(y))^{2\mu+1} (1-t^2)^\mu \dd t.
$
Thus, integration over each section reduces, up to a constant
factor, to integration on $[-1,1]$ with the Jacobi weight
$(1-t^2)^\mu$. This explains the use of the polynomials
$P_j^{(\mu,\mu)}$ on these sections. With $q_j$ denoting the orthonormal Jacobi polynomial introduced
at the beginning of this section (see \eqref{eq:qorth}), define
\begin{equation}\label{eq:Fj}
 F_j(s,y) = (a(y))^j  q_j\!\left(\frac{s}{a(y)}\right),
 \qquad \|y\|<1.
\end{equation}
Although the expression contains a square root, it extends to a
polynomial on $\R^d$.
Indeed, the parity relation $q_j (-t) = (-1)^j q_j (t)$ gives
\[
 q_j(t) = \sum_{r=0}^{\lfloor j/2 \rfloor} c_{j,r} t^{j-2r},
 \qquad c_{j,0} > 0,
\]
and consequently
\begin{equation}\label{eq:Fpolynomial}
 F_j(s,y) = \sum_{r=0}^{\lfloor j/2 \rfloor}
            c_{j,r} s^{j - 2r} (1 - \|y\|^2)^r.
\end{equation}
Each summand has total degree at most $j$, and the coefficient
of $s^j$ is $c_{j,0} \ne 0$. Hence $F_j$ has total degree $j$.
Formula \eqref{eq:Fpolynomial} defines $F_j$ on $\R^d$.

\begin{lemma}\label{lem:fiberdecomposition}
Every $p \in \PP_n^d$ has a unique representation
\begin{equation}\label{eq:fiberexpansion}
 p(s,y) = \sum_{j=0}^n F_j (s,y) g_j (y),
 \qquad g_j \in \PP_{n-j}^{d-1}.
\end{equation}
If $p, \widetilde p$ have coefficient polynomials
$ g_j, \widetilde g_j$, respectively, then
\begin{align}
 &\ip{p}{\widetilde p}
   =\sum_{j=0}^n \int_{\B^{d-1}}
      g_j(y) \widetilde g_j(y)
      (1 - \|y\|^2)^{\mu + j + 1/2} \dd y,
      \label{eq:fibermass}\\
 &\ip{\Phi_1\partial_sp}{\Phi_1\partial_s\widetilde p}
   =\sum_{j=1}^n \eta_j \int_{\B^{d-1}}
      g_j(y) \widetilde g_j(y)
      (1 - \|y\|^2)^{\mu + j - 1/2} \dd y.
      \label{eq:fiberenergy}
\end{align}
\end{lemma}

\begin{proof}
Recall that $x = (s,y)$, where $s = x_1$ and $y = (x_2,\ldots,x_d)$.
Let $p \in \PP_n^d$ be nonzero, and let $k$ denote its degree
with respect to the variable $s$. Then
\[
 p(s,y)=\sum_{j=0}^k s^j b_j(y),\qquad
 b_j \in \PP_{n-j}^{d-1}, \quad b_k \not \equiv0.
\]
By \eqref{eq:Fpolynomial}, the coefficient of $s^k$ in
$F_k(s,y)$ is $c_{k,0}$. Hence the polynomial
\[
 p_1(s,y)
 =p(s,y) - \frac{b_k(y)}{c_{k,0}} F_k(s,y)
\]
belongs to $\PP_n^d$ and is either zero or has degree
strictly less than $k$ with respect to $s$.
The total degree remains at most $n$, because
$\deg F_k = k$ and $\deg b_k \le n-k$.
Repeating the procedure with the remainder gives
\eqref{eq:fiberexpansion} and the stated bounds on $\deg g_j$.
To prove uniqueness, suppose that
\[
 p(s,y)=\sum_{j=0}^n F_j(s,y) g_j(y)
       =\sum_{j=0}^n F_j(s,y) \widehat g_j(y).
\]
If the two representations differ, let $k$ be the largest
index for which $g_k\ne\widehat g_k$. Then
\[
 \sum_{j=0}^k F_j(s,y) \left(g_j(y)-\widehat g_j(y)\right) = 0.
\]
For $j<k$, the degree of $F_j$ in $s$ is less than $k$,
whereas the coefficient of $s^k$ in $F_k$ is $c_{k,0}$.
Comparing the coefficients of $s^k$ therefore gives
$
 c_{k,0} \left( g_k(y) - \widehat g_k(y) \right)=0.
$
Since $c_{k,0} > 0$, this contradicts the choice of $k$.
Hence the representation is unique.

To calculate the integrals, fix $j,k$ and consider the terms
$F_j(s,y) g_j(y)$ and $F_k(s,y) \widetilde g_k(y)$.
For brevity, write $g = g_j$ and $\widetilde g = \widetilde g_k$.
For each fixed $y$ with $\|y\| < 1$,
\eqref{eq:sectionchange} gives
\begin{align*}
 &\int_{-a(y)}^{a(y)}
       F_j(s,y) F_k(s,y) g(y) \widetilde g(y)
       (1 - s^2 - \|y\|^2)^\mu\dd s\\
 &\qquad=
 g(y) \widetilde g(y) a(y)^{j+k+2\mu+1}
       \int_{-1}^1 q_j(t) q_k(t) (1-t^2)^\mu \dd t =
 \delta_{jk}\, g(y) \widetilde g(y) a(y)^{2j+2\mu+1}.
\end{align*}
Let
$
 p=\sum_{j=0}^n F_j g_j,
$
$
 \widetilde p = \sum_{k=0}^n F_k \widetilde g_k.
$
By Fubini's theorem and the preceding calculation,
\begin{align*}
 \int_{\B^d} p(x) \widetilde p(x) W_\mu(x) \dd x
 &=
 \sum_{j,k = 0}^n
 \int_{\B^{d-1}} g_j(y) \widetilde g_k(y)\\
 &\qquad\times
 \left[
 \int_{-a(y)}^{a(y)}
 F_j(s,y) F_k(s,y) (1 - s^2 - \|y\|^2)^\mu \dd s
 \right]\dd y\\
 &=
 \sum_{j,k=0}^n \delta_{jk}
 \int_{\B^{d-1}}
 g_j(y) \widetilde g_k(y) a(y)^{j+k+2\mu+1} \dd y.
\end{align*}
Consequently,
\[
 \int_{\B^d} p(x) \widetilde p(x) W_\mu(x) \dd x
 =
 \sum_{j=0}^n\int_{\B^{d-1}}
 g_j(y) \widetilde g_j(y) (1-\|y\|^2)^{\mu + j + \frac12} \dd y,
\]
which proves \eqref{eq:fibermass}.

For $j \ge 1$, differentiation with respect to $s$ gives
$
 \partial_s F_j(s,y) = a(y)^{j-1} q_j'(s/a(y)).
$
Thus for $j,k \ge 1$ the corresponding derivative integral on the
section is
\begin{align*}
 &\int_{-a(y)}^{a(y)}
       (\Phi_1(s,y))^2
       \partial_s F_j(s,y) \partial_s F_k(s,y)
       g(y) \widetilde g(y) (1-s^2-\|y\|^2)^\mu \dd s\\
 &\qquad=
 g(y) \widetilde g(y) (a(y))^{j+k+2\mu-1}
       \int_{-1}^1 q_j'(t) q_k'(t) (1-t^2)^{\mu+1} \dd t
 = \eta_j \delta_{jk}\, g(y) \widetilde g(y) (a(y))^{2j+2\mu-1}.
\end{align*}
The last equality follows from Lemma~\ref{lem:jacobi}.
Integrating over $\B^{d-1}$ and summing over $j,k$,
we obtain \eqref{eq:fiberenergy}.
\end{proof}
We next decompose each polynomial $g_j$ into harmonic and radial components.
For $d \ge 3$, let $\HH_\ell^{d-1}$ denote the space
of homogeneous harmonic polynomials of degree $\ell$ on $\R^{d-1}$, and set
\[
 h_\ell^{d-1} = \dim\HH_\ell^{d-1}
 =\binom{\ell+d-2}{d-2} - \binom{\ell+d-4}{d-2},
\]
where the second binomial coefficient is understood to be zero for $\ell=0,1$.
Choose a basis $Y_{\ell,\nu}$, $1 \le \nu \le h_\ell^{d-1}$,
whose restrictions to $\Sph^{d-2}$ are orthonormal for surface measure:
\begin{equation}\label{eq:sphericalorthogonality}
 \int_{\Sph^{d-2}} Y_{\ell,\nu}(\theta)
          Y_{k,\rho}(\theta) \dd\sigma(\theta)
       =\delta_{\ell k} \delta_{\nu\rho}.
\end{equation}
For $d \ge 3$, the restriction of $Y_\ell \in \HH_\ell^{d-1}$
to $\mathbb S^{d-2}$ satisfies
$
 -\Delta_{\mathbb S^{d-2}} Y_\ell
 =\ell(\ell + d - 3)Y_\ell;
$
see \cite[Chapter~1]{DX}.
For $d=2$, $h_0^1 = h_1^1 = 1$ and
\begin{equation}\label{eq:D1basis}
 Y_{0,1}(y)= \frac1{\sqrt2},\qquad
 Y_{1,1}(y)= \frac{y}{\sqrt2}.
\end{equation}

\begin{lemma}\label{lem:harmonic}
Let $d \ge 2$ and $m \ge 0$ be integers.
Every $g \in \PP_m^d$ has a unique representation
\begin{equation}\label{eq:harmonicexpansion}
 g(y)=
 \sum_{\ell=0}^{m}
 \sum_{\nu=1}^{h_\ell^d}
 Y_{\ell,\nu}(y)h_{\ell,\nu}(\|y\|^2),
\qquad
\text{where}
\qquad
h_{\ell,\nu} \in
\PP_{\lfloor(m-\ell)/2\rfloor}^1.
\end{equation}
The summands in \eqref{eq:harmonicexpansion} corresponding
to distinct pairs $(\ell,\nu)$ are orthogonal with respect
to every radial weight on $\B^d$ for which they are
square integrable.

Moreover, for $\gamma > -1$, $\ell=0,1,\ldots$,
$1 \le \nu,\rho \le h_\ell^d$, and univariate polynomials $h,k$,
\begin{align}
 &\int_{\B^d}
 Y_{\ell,\nu}(y) Y_{\ell,\rho}(y)
 h(\|y\|^2) k(\|y\|^2) (1 - \|y\|^2)^\gamma\dd y
 =
 \frac{\delta_{\nu\rho}}{2}
 \int_0^1 h(u) k(u)
 u^{\ell+d/2-1} (1-u)^\gamma\dd u.
 \label{eq:radialintegral}
\end{align}
\end{lemma}

\begin{proof}
The representation \eqref{eq:harmonicexpansion}, with the
stated degree bounds, is given in
\cite[Proposition~1]{Xu2004Interpolation}.
Its uniqueness follows from the uniqueness of the harmonic
decomposition of each homogeneous part of $g$ and the choice
of a basis in each space $\HH_\ell^d$; see
\cite[Section~1.1]{DX}.

Let $w(\|y\|)$ be a radial weight for which the polynomials
under consideration are square integrable.
For $y=r\theta$, where $0 < r < 1$ and
$\theta \in \Sph^{d-1}$, homogeneity gives
\[
 Y_{\ell,\nu}(r\theta)
 =r^\ell Y_{\ell,\nu}(\theta).
\]
Thus, integration in polar coordinates yields
\begin{align*}
 &\int_{\B^d}
 Y_{\ell,\nu}(y)  Y_{\ell',\rho}(y)
 h(\|y\|^2)  k(\|y\|^2)  w(\|y\|)  \dd y\\
 &\quad=
 \left(
 \int_{\Sph^{d-1}}
 Y_{\ell,\nu}(\theta)  Y_{\ell',\rho}(\theta)
 \dd\sigma(\theta)
 \right)
 \int_0^1
 h(r^2)  k(r^2)  r^{\ell+\ell'+d-1}  w(r)  \dd r.
\end{align*}
By orthonormality on the sphere, the first factor equals
$\delta_{\ell\ell'}\delta_{\nu\rho}$.
Hence the integral vanishes whenever
$(\ell,\nu)\ne(\ell',\rho)$, proving the orthogonality
assertion.

Finally, taking $\ell' = \ell$ and
$w(r) = (1 - r^2)^\gamma$, we obtain
\begin{align*}
 \int_{\B^d}
 Y_{\ell,\nu}(y)  Y_{\ell,\rho}(y)
 h(\|y\|^2)  k(\|y\|^2)  (1 - \|y\|^2)^\gamma  \dd y
 &=
 \delta_{\nu\rho}
 \int_0^1
 h(r^2)  k(r^2)  r^{2\ell+d-1}  (1-r^2)^\gamma  \dd r\\
 &=
 \frac{\delta_{\nu\rho}}{2}
 \int_0^1
 h(u)  k(u)  u^{\ell+d/2-1} (1-u)^\gamma \dd u.
\end{align*}
This proves \eqref{eq:radialintegral}.
\end{proof}

In dimension one, the corresponding decomposition is obtained
by separating the even and odd parts. Every $g\in\PP_m^1$
can be written uniquely as
\[
 g(y) = Y_{0,1}(y)  h_0(y^2) + Y_{1,1}(y)  h_1(y^2),
\]
where $Y_{0,1}$ and $Y_{1,1}$ are given by
\eqref{eq:D1basis},
$\deg h_0 \le \lfloor m/2\rfloor$, and, for $m \ge 1$,
$\deg h_1 \le \lfloor(m-1)/2\rfloor$.
For $m=0$, we set $h_1=0$.
The two summands are orthogonal with respect to every
even weight for which they are square integrable.
Moreover, for $\ell \in \{0,1\}$, the substitution $u = y^2$
gives \eqref{eq:radialintegral} with $d=1$.

For $(j,\ell) \in \mathcal I_{n,d}$ and
$1 \le \nu \le h_\ell^{d-1}$, consider polynomials
\begin{equation}\label{eq:blockpolynomial}
 p_h(s,y) = F_j(s,y)  Y_{\ell,\nu}(y)  h(\|y\|^2),
 \qquad h \in \PP_{m_{j,\ell}}^1.
\end{equation}
Their total degree is at most
$j + \ell + 2m_{j,\ell} \le n$.
By Lemmas~\ref{lem:fiberdecomposition} and~\ref{lem:harmonic},
\begin{align}
 &\ip{p_h}{p_k}
  =\frac12 \int_0^1
      h(u)  k(u)  u^{\beta_\ell}  (1-u)^{\alpha_j+1}  \dd u,
      \label{eq:blockmass}\\
 &\ip{\Phi_1\partial_sp_h}{\Phi_1\partial_sp_k}
  =\frac{\eta_j}2  \int_0^1
      h(u)  k(u)  u^{\beta_\ell}  (1-u)^{\alpha_j}  \dd u.
      \label{eq:blockenergy}
\end{align}
Here
$\ell + (d-1)/2 - 1 = \beta_\ell$,
$\mu + j + 1/2 = \alpha_j + 1$, and
$\mu + j - 1/2 = \alpha_j$.
Polynomials corresponding to different triples
$(j,\ell,\nu)$ are orthogonal for both integrals.
This includes orthogonality of all $j \ge 1$ terms to the
$j=0$ part of \eqref{eq:fiberexpansion} for $\ip{\cdot}{\cdot}$.
The latter part consists of polynomials independent of $s$.

For a nonzero polynomial of the form \eqref{eq:blockpolynomial},
we have
\begin{equation}\label{eq:blockquotient}
 \frac{\nm{\Phi_1\partial_sp_h}^{\,2}}{\nm{p_h}^{\,2}}
 =
 \eta_j
 \frac{\int_0^1 (h(u))^2 u^{\beta_\ell} (1-u)^{\alpha_j} \dd u}
      {\int_0^1 (h(u))^2 u^{\beta_\ell} (1-u)^{\alpha_j+1} \dd u}.
\end{equation}

\section{The sharp constant and all equality cases}
\label{sec:sharp}

For each $(j,\ell)\in\mathcal I_{n,d}$, let
$
 R_{j,\ell}(u)=
 P_{N_{j,\ell}}^{(\alpha_j,\beta_\ell)}(2u-1),
$
$
 u_{j,\ell,q}=\frac{1+\xi_{j,\ell,q}}2.
$
Define
\begin{align}
 &L_{j,\ell,q}(u)
  =\frac{R_{j,\ell}(u)}
        {(u-u_{j,\ell,q}) R_{j,\ell}'(u_{j,\ell,q})},
       \label{eq:indexedcardinal}\\
 &\gamma_{j,\ell,q}
  =\int_0^1 (L_{j,\ell,q}(u))^2
          u^{\beta_\ell} (1-u)^{\alpha_j} \dd u.
       \label{eq:indexedgamma}
\end{align}
Let
$
 \mathcal K_n = \{p \in \PP_n^d: \partial_sp=0\}
              = \{p(s,y)=g(y): g \in \PP_n^{d-1}\}.
$

Using the polynomials $F_j$ defined in \eqref{eq:Fj}
and the fixed basis
$\{Y_{\ell,\nu}\}_{\nu=1}^{h_\ell^{d-1}}$ of $\HH_\ell^{d-1}$,
normalized as in \eqref{eq:sphericalorthogonality},
we obtain the following polynomial expansion, which gives
explicit formulas for the numerator and denominator
in \eqref{eq:Cdef}.

\begin{theorem}\label{thm:commonbasis}
For $(j,\ell)\in\mathcal I_{n,d}$,
$1 \le \nu \le h_\ell^{d-1}$, and $1 \le q \le N_{j,\ell} = m_{j,\ell}+1$, let
\begin{equation}\label{eq:commonbasis}
 \psi_{j,\ell,\nu,q}(s,y)=
 \left(\frac{2}
 {\gamma_{j,\ell,q} (1-u_{j,\ell,q})} \right)^{1/2}
 F_j(s,y) Y_{\ell,\nu}(y) L_{j,\ell,q}(\|y\|^2).
\end{equation}
These polynomials are orthonormal with respect to $\ip{\cdot}{\cdot}$
and are orthogonal to $\mathcal K_n$.
Every $p \in \PP_n^d$ has a unique expansion
\begin{equation}\label{eq:fullpolynomialexpansion}
 p = p_0+
   \sum_{(j,\ell) \in \mathcal I_{n,d}}
   \sum_{\nu=1}^{h_\ell^{d-1}} \sum_{q=1}^{N_{j,\ell}}
       c_{j,\ell,\nu,q} \psi_{j,\ell,\nu,q},
 \qquad p_0 \in \mathcal K_n.
\end{equation}
For this expansion,
\begin{align}
 &\nm{p}^{\,2}
  =\nm{p_0}^{\,2}
     +\sum_{j,\ell,\nu,q} c_{j,\ell,\nu,q}^2,
       \label{eq:massdiagonal}\\
 &\nm{\Phi_1\partial_sp}^{\,2}
  =\sum_{j,\ell,\nu,q}
       \Lambda_{j,\ell,q}c_{j,\ell,\nu,q}^2,
       \label{eq:energy-diagonal}
\end{align}
where the sums are taken over the same index ranges as in
\eqref{eq:fullpolynomialexpansion}, and
\begin{equation}\label{eq:quotientvalues}
 \Lambda_{j,\ell,q}=
 \frac{\eta_j}{1-u_{j,\ell,q}}
 =\frac{2j(j+2\mu+1)}{1-\xi_{j,\ell,q}}.
\end{equation}
\end{theorem}

\begin{proof}
By \eqref{eq:Fpolynomial}, $F_j$ has total degree $j$.
Since $Y_{\ell,\nu}$ is homogeneous of degree $\ell$
and $L_{j,\ell,q}$ has degree $m_{j,\ell}$, we obtain
\[
 \deg\psi_{j,\ell,\nu,q}
 =j + \ell + 2m_{j,\ell} \le n.
\]
Thus $\psi_{j,\ell,\nu,q} \in \PP_n^d$.
Fix an admissible triple $(j,\ell,\nu)$ and write
$L_q=L_{j,\ell,q}$, $u_q=u_{j,\ell,q}$, and
$\gamma_q=\gamma_{j,\ell,q}$.
Let
\[
 p_q(s,y)=F_j(s,y)Y_{\ell,\nu}(y)L_q(\|y\|^2).
\]
Equations \eqref{eq:blockmass}, \eqref{eq:blockenergy},
and \eqref{eq:cardinalorthogonality} give
\begin{align*}
 &\ip{p_q}{p_r}
     =\tfrac12 \gamma_q(1-u_q) \delta_{qr},\\
 &\ip{\Phi_1\partial_sp_q}{\Phi_1\partial_sp_r}
     =\tfrac12 \eta_j \gamma_q \delta_{qr}.
\end{align*}
The normalization in \eqref{eq:commonbasis} therefore gives
unit norm and derivative norm squared
$\eta_j/(1-u_q)=\Lambda_{j,\ell,q}$.
For different triples $j,\ell,\nu$, both mixed integrals
are zero by Section~\ref{sec:lemmas}.
It follows that, for any two admissible quadruples of indices,
\begin{align*}
 &\ip{\psi_{j,\ell,\nu,q}}{\psi_{j',\ell',\nu',q'}}
     =\delta_{jj'}\delta_{\ell\ell'}
       \delta_{\nu\nu'}\delta_{qq'},\\
 &\ip{\Phi_1\partial_s\psi_{j,\ell,\nu,q}}
       {\Phi_1\partial_s\psi_{j',\ell',\nu',q'}}
     =\Lambda_{j,\ell,q}
       \delta_{jj'}\delta_{\ell\ell'}
       \delta_{\nu\nu'}\delta_{qq'}.
\end{align*}

To prove completeness, start with the unique expansion
\eqref{eq:fiberexpansion}.
Its $j=0$ term is $p_0 = F_0 g_0 \in \mathcal K_n$;
conversely, all elements of $\mathcal K_n$ are of this form,
because $F_0 = q_0$ is a nonzero constant.
For each $j\ge1$, apply Lemma~\ref{lem:harmonic} in
dimension $d-1$ when $d\ge3$, and use the even--odd
decomposition following that lemma when $d=2$.
Every resulting one-variable coefficient polynomial has
degree at most $m_{j,\ell}$, and hence has a unique expansion
in the basis $L_{j,\ell,q}$.
This yields \eqref{eq:fullpolynomialexpansion}.
Uniqueness follows successively from the two decomposition
lemmas and from uniqueness of interpolation at the distinct
nodes $u_{j,\ell,q}$.
Since $p_0$ is independent of $s$, it is orthogonal to
all $\psi_{j,\ell,\nu,q}$ by \eqref{eq:fibermass}.
Using the above orthogonality relations and noting that
$\partial_s p_0 = 0$, we obtain \eqref{eq:massdiagonal}
and \eqref{eq:energy-diagonal}.
\end{proof}

\begin{proof}[Proof of Theorem~\ref{thm:main}]
It suffices to consider $i=1$. Let
\[
 M=\max_{(j,\ell) \in \mathcal I_{n,d}}
       \Lambda_{j,\ell,N_{j,\ell}}.
\]
For each fixed pair $(j,\ell)$, the zeros
$u_{j,\ell,q}$ increase strictly with $q$.
Since $\eta_j > 0$, formula \eqref{eq:quotientvalues} gives
\[
 \Lambda_{j,\ell,q}
 =\frac{\eta_j}{1-u_{j,\ell,q}}
 \le \Lambda_{j,\ell,N_{j,\ell}}
 \le M.
\]
Consequently, \eqref{eq:massdiagonal} and
\eqref{eq:energy-diagonal} imply that
\[
 \nm{\Phi_1\partial_s p}^{\,2}
 =\sum_{j,\ell,\nu,q}
       \Lambda_{j,\ell,q} c_{j,\ell,\nu,q}^2
 \le M \sum_{j,\ell,\nu,q} c_{j,\ell,\nu,q}^2
 \le M \nm{p}^{\,2}
\]
for every $p \in \PP_n^d$.
By \eqref{eq:Cdef}, this proves $C_{n,d,\mu}^2 \le M$.
Conversely, choose $(j_0,\ell_0) \in \mathcal I_{n,d}$
such that
$
 \Lambda_{j_0,\ell_0,N_{j_0,\ell_0}} = M,
$
and let
$
 p_* = \psi_{j_0,\ell_0,1,N_{j_0,\ell_0}}.
$
By Theorem~\ref{thm:commonbasis},
\[
 \nm{p_*}^{2} = 1,
 \qquad
 \nm{\Phi_1\partial_s p_*}^{2} = M.
\]
Thus the quotient in \eqref{eq:Cdef} equals $M$ for
$p = p_*$, so $C_{n,d,\mu}^2 \ge M$.
Hence $C_{n,d,\mu}^2 = M$.
The result for the other coordinate derivatives follows by permuting coordinates.
\end{proof}

\begin{corollary}\label{cor:extremizers}
Let
$
 \mathcal M_{n,d,\mu}
 =
 \left\{
 (j,\ell) \in \mathcal I_{n,d}:
 \Lambda_{j,\ell,N_{j,\ell}} = C_{n,d,\mu}^2
 \right\}.
$ For $i = 1$, a nonzero polynomial $p \in \PP_n^d$
attains the maximum in \eqref{eq:Cdef} if and only if
it has the representation
\begin{equation}\label{eq:all-extremizers}
 p(s,y)=
 \sum_{(j,\ell) \in \mathcal M_{n,d,\mu}}
 F_j(s,y) H_{j,\ell}(y)
 \frac{R_{j,\ell}(\|y\|^2)}
      {\|y\|^2-u_{j,\ell,N_{j,\ell}}},
\end{equation}
where $H_{j,\ell} \in \HH_\ell^{d-1}$ and at least one
of these polynomials is nonzero.

For $2\le i\le d$, the extremal polynomials are obtained
from those in \eqref{eq:all-extremizers} by interchanging
the coordinates $x_1$ and $x_i$.
\end{corollary}

\begin{remark}\label{rem:multiplicities}
The identities in Theorem~\ref{thm:commonbasis} retain every
coefficient of the derivative integral, not only its maximum.
For a fixed $j,\ell,q$, the coefficient
$\Lambda_{j,\ell,q}$ occurs $h_\ell^{d-1}$ times,
once for each harmonic basis element.
If different indices give the same value, the numbers of
corresponding basis polynomials add.
There are
\[
 \sum_{(j,\ell) \in \mathcal I_{n,d}}
 h_\ell^{d-1}(m_{j,\ell}+1)
 =\sum_{j=1}^n \dim\PP_{n-j}^{d-1}
 =\binom{n+d-1}{d}
\]
such polynomials.
Together with
$\dim\mathcal K_n = \binom{n+d-1}{d-1}$,
this accounts for all of $\PP_n^d$.
\end{remark}

We next compare the sharp constant with the upper bound obtained
from the classical Bernstein inequality \eqref{eq:classicalball}.
For polynomials $p,q$, define
\begin{equation}\label{eq:energyform}
 \mathcal E_\mu(p,q)=
 \int_{\B^d}
 \bigl(\nabla p \cdot \nabla q
       -(x\cdot\nabla p) (x\cdot\nabla q)\bigr) W_\mu(x) \dd x.
\end{equation}
The identity
\begin{align*}
 \sum_{1 \le i < j \le d} D_{ij}p\,D_{ij}q
 &=
 \frac12\sum_{i,j=1}^d
 (x_i \partial_j p-x_j \partial_ip)
 (x_i \partial_j q-x_j \partial_iq)\\
 &=
 \sum_{i,j=1}^d x_i^2\,\partial_jp\,\partial_jq
 -
 \sum_{i,j=1}^d x_ix_j\,\partial_ip\,\partial_jq\\
 &=
 \left(\sum_{i=1}^d x_i^2 \right)
 \left(\sum_{j=1}^d \partial_jp\,\partial_jq\right)
 -
 \left(\sum_{i=1}^d x_i \partial_ip\right)
 \left(\sum_{j=1}^d x_j \partial_jq\right)\\
 &=
 \|x\|^2\,\nabla p \cdot \nabla q
 -(x \cdot \nabla p)(x \cdot \nabla q)
\end{align*}
shows that
\[
 \mathcal E_\mu(p,p)=
 \sum_{i=1}^d \nm{\sqrt{1-\|x\|^2}\,\partial_ip}^{2}
 +\sum_{1 \le i < j \le d}\nm{D_{ij}p}^{2}.
\]
In the above notation, inequality~\eqref{eq:classicalball} takes the form
\begin{equation}\label{eq:energybound}
 \mathcal E_\mu(p,p)
 \le n(n + 2 \mu + d)\nm{p}^{2}.
\end{equation}
The factors $\Phi_i$ are closely related to the geometric
Bernstein inequalities of Baran~\cite{B1,B2}.
The following lemma relates these factors to the quadratic
form defining $\mathcal E_\mu$.
\begin{lemma}\label{lem:energybound}
For $p \in \PP_n^d$ and $\|x\| < 1$,
\begin{equation}\label{eq:pointwisecomparison}
 (\Phi_i(x))^2 |\partial_i p(x)|^2
 \le \|\nabla p(x)\|^2 - |x \cdot \nabla p(x)|^2.
\end{equation}
\end{lemma}

\begin{proof}
Fix $x \in \mathbb R^d$ with $\|x\|<1$ and define
$
 A = I-xx^{\mathsf T}.
$
For every $v \in \R^d$,
\[
 v^{\mathsf T} Av
 =\|v\|^2 - (x\cdot v)^2
 \ge (1 - \|x\|^2) \|v\|^2,
\]
so $A$ is positive definite. The Sherman–Morrison formula yields
\[
 A^{-1} = I + \frac{xx^{\mathsf T}}{1-\|x\|^2}.
\]
Hence, denoting by $e_i$ the $i$-th coordinate vector,
we obtain
\[
 e_i^{\mathsf T}A^{-1}e_i
 =1 + \frac{x_i^2}{1-\|x\|^2}
 =\frac{1 - \|x\|^2 + x_i^2}{1 - \|x\|^2}
 =(\Phi_i(x))^{-2}.
\]
Since $A^{-1/2}$ is symmetric, we have
\begin{align*}
 \bigl(A^{-1/2}e_i\bigr) \cdot
 \bigl(A^{1/2}\nabla p(x)\bigr)
 =
 e_i^{\mathsf T}(A^{-1/2})^{\mathsf T}
 A^{1/2} \nabla p(x)
 =
 e_i^{\mathsf T} A^{-1/2} A^{1/2} \nabla p(x)
 =\partial_i p(x).
\end{align*}
Consequently, the Cauchy--Schwarz inequality gives
\begin{align*}
 |\partial_i p(x)|^2
 \le
 \|A^{-1/2} e_i\|^2
 \|A^{1/2} \nabla p(x)\|^2
 &=
 \bigl(e_i^{\mathsf T} A^{-1} e_i\bigr)
 \bigl(\nabla p(x)^{\mathsf T} A \nabla p(x)\bigr)\\
 &=
 (\Phi_i(x))^{-2}
 \bigl(\|\nabla p(x)\|^2 - |x \cdot \nabla p(x)|^2\bigr).
\end{align*}
Multiplying by $(\Phi_i(x))^2$ proves
\eqref{eq:pointwisecomparison}.
\end{proof}

\begin{proposition}\label{prop:bounds}
For $d \ge 2$, $\mu > -1$, and $n \ge 1$,
\begin{equation}\label{eq:twosided}
 n(n+2\mu+1) \frac{n+\mu+d/2}{n+\mu+1/2}
 \le C_{n,d,\mu}^2 \le n(n + 2 \mu + d).
\end{equation}
In particular, for fixed $d,\mu$,
\begin{equation}\label{eq:asymptotic}
 \lim_{n \to \infty}\frac{C_{n,d,\mu}}n=1.
\end{equation}
\end{proposition}

\begin{proof}
Integrating \eqref{eq:pointwisecomparison} and using
\eqref{eq:energybound} gives the upper bound.
For the lower bound, take $p(s,y) = F_n(s,y)$.
By \eqref{eq:Fpolynomial}, this is a nonzero polynomial
of total degree $n$. In the expansion
\eqref{eq:fiberexpansion}, this choice corresponds to
\[
 g_n(y) \equiv 1,\qquad
 g_j(y) \equiv 0\quad (0 \le j < n).
\]
Since $g_n$ is constant, its harmonic decomposition
contains only the component of degree $\ell=0$,
with a constant radial factor.
Thus the relevant indices are $j=n$ and $\ell=0$, and
$
 m_{n,0}=0.
$
Therefore, \eqref{eq:blockquotient} gives
\[
 \frac{\nm{\Phi_1\partial_s F_n}^{2}}
      {\nm{F_n}^{2}}
 =
 \eta_n
 \frac{\displaystyle
       \int_0^1 u^{(d-3)/2} (1-u)^{n+\mu-1/2} \dd u}
      {\displaystyle
       \int_0^1 u^{(d-3)/2} (1-u)^{n+\mu+1/2} \dd u}.
\]
Recall the beta integral
\[
 B(a,b) = \int_0^1 u^{a-1} (1-u)^{b-1} \dd u
       =\frac{\Gamma(a) \Gamma(b)}{\Gamma(a+b)},
 \qquad a,b > 0.
\]
The elementary identity
$B(a,b+1) = bB(a,b)/(a+b)$ therefore gives
\begin{align*}
 \frac{\nm{\Phi_1\partial_sF_n}^{2}}{\nm{F_n}^{2}}
 =\eta_n\,
       \frac{B((d-1)/2,\mu+n+1/2)}
            {B((d-1)/2,\mu+n+3/2)}
 =n(n+2\mu+1)
       \frac{n+\mu+d/2}{n+\mu+1/2}.
\end{align*}
The asymptotic relation \eqref{eq:asymptotic} follows from these bounds.
\end{proof}

\section*{Acknowledgments}
The author was supported by the Polish National Science Centre (NCN) Miniatura grant no. 2025/09/X/ST1/00082.

The author used ChatGPT to assist with language editing and
to review the mathematical exposition. The author independently
checked the final manuscript and takes responsibility for its content.


\begin{thebibliography}{9}

\bibitem{B1}
M.~Baran,
\emph{Bernstein type theorems for compact sets in $\mathbb{R}^n$},
J. Approx. Theory \textbf{69} (1992), 156--166.
\doi{10.1016/0021-9045(92)90139-F}.

\bibitem{B2}
M.~Baran,
\emph{Bernstein type theorems for compact sets in $\mathbb{R}^n$
revisited},
J. Approx. Theory \textbf{79} (1994), 190--198.
\doi{10.1006/jath.1994.1124}.

\bibitem{TB}
T.~Beberok,
\emph{Weighted Bernstein-type inequalities on a trapezoidal domain},
preprint, 2026.
\href{https://arxiv.org/abs/2609.33283}{arXiv:2609.33283}.

\bibitem{BX}
T.~Beberok and Y.~Xu,
\emph{On Bernstein inequalities on the unit ball},
Constr. Approx. (2026), published online 10 June 2026.
\doi{10.1007/s00365-026-09761-8}.

\bibitem{Dai2006}
F.~Dai,
\emph{Multivariate polynomial inequalities with respect to
doubling weights and $A_\infty$ weights},
J. Funct. Anal. \textbf{235} (2006), 137--170.
\doi{10.1016/j.jfa.2005.09.009}.

\bibitem{DX2010}
F.~Dai and Y.~Xu,
\emph{Moduli of smoothness and approximation on the unit sphere
and the unit ball},
Adv. Math. \textbf{224} (2010), 1233--1310.
\doi{10.1016/j.aim.2010.01.001}.

\bibitem{DX}
F.~Dai and Y.~Xu,
\emph{Approximation Theory and Harmonic Analysis on Spheres and Balls},
Springer Monographs in Mathematics,
Springer, New York, 2013.
\doi{10.1007/978-1-4614-6660-4}.

\bibitem{Gautschi}
W.~Gautschi,
\emph{Orthogonal Polynomials: Computation and Approximation},
Numerical Mathematics and Scientific Computation,
Oxford University Press, Oxford, 2004.
\doi{10.1093/oso/9780198506720.001.0001}.

\bibitem{GX}
Y.~Ge and Y.~Xu,
\emph{Sharp Bernstein inequalities on simplex},
Constr. Approx. \textbf{62} (2025), 305--328.
\doi{10.1007/s00365-024-09680-6}.

\bibitem{GM}
A.~Guessab and G.~V.~Milovanovi\'c,
\emph{Weighted $L^2$-analogues of Bernstein's inequality
and classical orthogonal polynomials},
J. Math. Anal. Appl. \textbf{182} (1994), 244--249.
\doi{10.1006/jmaa.1994.1078}.

\bibitem{Kroo}
A.~Kro\'o,
\emph{Exact $L_2$ Bernstein--Markov inequality on the ball},
J. Approx. Theory \textbf{281--282} (2022), article~105795.
\doi{10.1016/j.jat.2022.105795}.

\bibitem{Kroo2}
A.~Kro\'o,
\emph{On multivariate $L^p$ Bernstein--Markov type inequalities},
Proc. Amer. Math. Soc. \textbf{152} (2024), 5149--5162.
\doi{10.1090/proc/16925}.

\bibitem{NS}
G.~Nikolov and A.~Shadrin,
\emph{Markov-type inequalities and extreme zeros of orthogonal
polynomials},
J. Approx. Theory \textbf{271} (2021), article~105644.
\doi{10.1016/j.jat.2021.105644}.

\bibitem{Szego}
G.~Szeg\H{o},
\emph{Orthogonal Polynomials},
4th ed., Amer. Math. Soc. Colloq. Publ., vol.~23,
American Mathematical Society, Providence, RI, 1975.
\doi{10.1090/coll/023}.

\bibitem{Xu2004Interpolation}
Y.~Xu,
\emph{On polynomial interpolation on the unit ball},
in: D.~K.~Dimitrov, G.~Nikolov, and R.~Uluchev (eds.),
\emph{Approximation Theory: A Volume Dedicated to Borislav Bojanov},
Marin Drinov Academic Publishing House, Sofia, 2004,
pp.~325--337.
\href{https://pages.uoregon.edu/yuan/paper/ball-inter.pdf}{Full text}.

\bibitem{Xu2023}
Y.~Xu,
\emph{Bernstein inequality on conic domains and triangles},
J. Approx. Theory \textbf{290} (2023), article~105889.
\doi{10.1016/j.jat.2023.105889}.

\end{thebibliography}
\end{document}